\documentclass[12pt,reqno]{amsart}
\usepackage{graphicx}
\usepackage{amsmath,amssymb,mathrsfs}
\usepackage[usenames]{color}

\newtheorem{Theorem}{Theorem}[section]

\newtheorem{Definition}[Theorem]{Definition}

\numberwithin{equation}{section}

\begin{document}

\title[]{ Schiffer's Conjecture, Interior Transmission Eigenvalues and Invisibility Cloaking: Singular Problem vs. Nonsingular Problem}
\author{Hongyu Liu}
\address{Department of
Mathematics and Statistics, University of North Carolina, Charlotte,
NC 28223, USA.}
\email{hongyu.liuip@gmail.com}

\maketitle

\begin{abstract}
In this note, we present some connections between Schiffer's conjecture, interior transmission eigenvalue problem and singular and non-singular invisibility cloaking problems of acoustic waves.  
\end{abstract}

\section{Schiffer's conjecture}

Schiffer's conjecture is a long standing problem in spectral theory, which is stated as follows:
\begin{quote}
Let $\Omega\subset\mathbb{R}^2$ be a bounded domain. Does the existence of a nontrivial solution $u$ to the over-determined Neumann eigenvalue problem
\begin{equation}\label{eq:Schiffer}
\begin{cases}
&-\Delta u=\lambda u\qquad\ \ \mbox{in\ \ $\Omega$}, \quad \lambda\in\mathbb{R}_+,\\
&\displaystyle{\frac{\partial u}{\partial n}}=0\hspace*{1.68cm} \mbox{on\ \ $\partial\Omega$},\\
& u=\mbox{const} \qquad\quad\,  \mbox{on\  \ $\partial\Omega$},
\end{cases}
\end{equation}
imply that $\Omega$ must be a ball?
\end{quote}

The problem is equivalent to the Pompeiu's problem in integral geometry. A domain $\Omega\subset\mathbb{R}^2$ is said to have the Pompeiu property if{f} the only continuous function $\varphi$ on $\mathbb{R}^2$ for which $\displaystyle{\int_{\sigma(\Omega)}\varphi(x,y)\ dxdy}=0$ for every rigid motion $\sigma$ of $\mathbb{R}^2$ is $\varphi(x,y)=0$.  
It is shown in \cite{Bro} that the failure of the Pompeiu property of a domain $\Omega$ is equivalent to the existence of a nontrivial solution to (\ref{eq:Schiffer}). We refer to \cite{Zal} for an extensive survey on the current state of the problem. 

For our subsequent discussion, we introduce the following theorem on the transformation invariance of the Laplace's equation, whose proof could be found, e.g., in \cite{GKLU5,KOVW,Liu}.

\begin{Theorem}\label{thm:transoptics}
Let $\Omega$ and $\widetilde{\Omega}$ be two bounded Lipschitz domains in $\mathbb{R}^N$ and suppose that there exists a diffeomorphism $F:\Omega\rightarrow \widetilde{\Omega}$. Let $u\in H^1(\Omega)$ satisfy 
\[
\nabla\cdot(g(x)\nabla u(x))+\lambda q(x) u(x)=f(x)\quad x\in \Omega,
\]
where $g(x)=(g^{ij}(x))_{i,j=1}^2, q(x)$, $x\in\Omega$ are uniformly elliptic and $g$ is symmetric, and $f\in L^2(\Omega)$. Then one has that $\widetilde{u}=(F^{-1})^* u:=u\circ F^{-1}\in H^1(\widetilde{\Omega})$ satisfies
\[
\nabla\cdot(\widetilde{g}(y)\nabla \widetilde{u}(y))+\lambda \widetilde{q}(y) \widetilde{u}(y)=\widetilde{f}(y), \ y\in \widetilde{\Omega},
\]
where
\begin{equation}\label{eq:trans}
\begin{split}
& F_*g(y)=\frac{DF(x)\cdot g(x)\cdot (DF(x))^T}{\left|\mbox{\emph{det}}(DF(x))\right|}\bigg|_{x=F^{-1}(y)},\\
& F_*q(y)=\frac{q(x)}{\left|\mbox{\emph{det}}(DF(x))\right|}\bigg|_{x=F^{-1}(y)},\ \ x\in\Omega, \ y\in\widetilde{\Omega},
\end{split}
\end{equation}
and
\[
\widetilde{f}=\left(\frac{f}{|\det(DF)|}\right)\circ F^{-1},
\]
and $DF$ denotes the Jacobian matrix of the transformation $F$. 
\end{Theorem}

Next, let $u$ and $\Omega$ be the ones in (\ref{eq:Schiffer}) and suppose there exists a diffeomorphism
\[
F:\Omega\rightarrow \widetilde{\Omega}.
\]
Further, we let
\begin{equation}\label{eq:trans2}
\widetilde{g}=F_* I\qquad \mbox{and}\qquad \widetilde{q}=F_* 1.
\end{equation}
Then, according to Theorem~\ref{thm:transoptics}, if one lets $v=(F^{-1})^* u$, we have
\begin{equation}\label{eq:general schiffer}
\begin{cases}
&\displaystyle{\nabla\cdot(\widetilde{g}\nabla v)}+\lambda \widetilde{q} v=0\quad\ \mbox{in\ \ $\widetilde{\Omega}$},\\
&\displaystyle{\sum_{i,j=1}^2 n_i\widetilde{g}^{ij}\partial_j v}=0\hspace*{1.2cm}\mbox{on\ \ $\partial\widetilde{\Omega}$},\\
&\ v=\mbox{const}\hspace*{2.3cm}\mbox{on\ \ $\partial\widetilde{\Omega}$}.
\end{cases}
\end{equation}
Hence, it is natural to propose the following generalized Schiffer's conjecture:
\begin{quote}
Let $(\Omega; g, q)$ be uniformly elliptic and symmetric.  Does the existence of a nontrivial solution to
\begin{equation}\label{eq:general schiffer 2}
\begin{cases}
&\displaystyle{\nabla\cdot({g}\nabla u)}+\lambda {q} u=0\quad\ \mbox{in\ \ ${\Omega}$},\ \ \lambda\in\mathbb{R}_+,\\
&\displaystyle{\sum_{i,j=1}^2 n_i {g}^{ij}\partial_j u}=0\hspace*{1.2cm}\mbox{on\ \ $\partial{\Omega}$},\\
&\ u=\mbox{const}\hspace*{2.3cm}\mbox{on\ \ $\partial\Omega$}.
\end{cases}
\end{equation}
imply that there must exist a diffeomorphism $F$ such that
\[
F: B\rightarrow \Omega,\quad \mbox{$B$ is a ball},
\]
and
\[
g=F_* I\qquad\mbox{and}\qquad q=F_*1?
\]
\end{quote} 

\begin{Definition}\label{def:1}
Let $(\Omega; g, q)$ be uniformly elliptic. Then, it is said to possess the Pompeiu property if the corresponding over-determined system (\ref{eq:general schiffer 2})  has only the trivial solution. 
\end{Definition}
Later in Section 3, we shall see that the Pompeiu property of the parameters $(\Omega; g,p)$ would have important implications for invisibility cloaking in acoustics.  

\section{Interior transmission eigenvalue problem}

Let $(\Omega; g_1, q_1)$ and $(\Omega; g_2, q_2)$ be uniformly elliptic and symmetric such that
\[
(g_1, q_1)\neq (g_2, q_2).
\]
Consider the following interior transmission problem for a pair $(u,v)$
\begin{equation}\label{eq:TE}
\begin{cases}
& \nabla\cdot(g_1\nabla u)+\lambda q_1 u=0\qquad\hspace*{2.5cm}\ \mbox{in\ $\Omega$},\\
& \nabla\cdot(g_2\nabla v)+\lambda q_2 v=0\qquad\hspace*{2.5cm}\ \,\mbox{in\ $\Omega$},\\
& a_{11} u+a_{22}v=0\hspace*{4.6cm}\,\mbox{on\  $\partial\Omega$},\\
& \displaystyle{a_{21}\sum_{i,j=1}^2 n_i g_1^{ij}\partial_j u+a_{22}\sum_{i,j=1}^2 n_i g_2^{ij}\partial_j v=0\quad\,\,\mbox{on\  $\partial\Omega$}},
\end{cases}
\end{equation}
where $A=(a_{ij})_{i,j=1}^2\in C(\partial\Omega)$.
If there exists a nontrivial pair of solutions $(u,v)$ to (\ref{eq:TE}), then $\lambda$ is called a generalized interior transmission eigenvalue, and $(u,v)$ are called generalized interior transmission eigenfunctions.

The interior transmission eigenvalue problem arises in inverse scattering theory and has a long history in literature (cf. \cite{ColKre,ColPaiSyl}). It was first introduced in \cite{ColMon} in connection with an inverse scattering problem for the reduced wave equation. Later, it found important applications in inverse scattering theory, especially, for the study of qualitative reconstruction schemes including linear sampling method and factorization method (see, e.g., \cite{Kir}). 

We would like to emphasize that the one (\ref{eq:TE}) presented here is a generalized formulation of the interior transmission eigenvalue problem that has been considered in the literature. 

\section{Transformation optics and invisibility cloaking}

In recent years, transformation optics and invisibility cloaking have received significant attentions;  see, e.g., \cite{CC,GKLU4,GKLU5,Nor} and references therein.  The crucial ingredient is the transformation invariance of the acoustic wave equations, which is actually Theorem~\ref{thm:transoptics}. Let's consider the well-known two-dimensional cloaking problem for the time-harmonic acoustic wave governed by the Helmholtz equation (cf. \cite{GLU2,Leo,PenSchSmi})
\begin{equation}\label{eq:Helm1}
\nabla\cdot(g\nabla u)+\omega^2 q u=f\chi_{B_1}\quad \mbox{in\ $B_2$},
\end{equation}
where $\omega\in\mathbb{R}_+$ is the wave frequency, $g$ is symmetric matrix-valued denoting the acoustical density and $q$ is the bulk modulus. In (\ref{eq:Helm1}), the parameters are given as
\begin{equation}\label{eq:parameter}
(g, q)=\begin{cases}
&g_a,\ q_a\quad\mbox{in\ $B_1$},\\
& g_c,\ q_c\quad\mbox{\ in\ $B_2\backslash B_1$},
\end{cases}
\end{equation}
and $f\in L^2(B_1)$. In the physical situation, $(g_c, q_c)$ is the cloaking medium which we shall specify later, and $(g_a, q_a, f_a)$ is the target object including a passive medium $(g_a, q_a)$ and an active source or sink $f_a$. The invisibility is understood in terms of the exterior measurement encoded either into the boundary Dirichlet-to-Neumann (DtN) map or the scattering amplitude far-away from the object. To ease our exposition, we only consider the boundary DtN map. Let 
\begin{equation}\label{eq:bc}
u=\varphi\in H^{1/2}(\partial B_2)
\end{equation}
and define the DtN operator by
\[
\Lambda_c(\varphi)=\sum_{i,j=1}^2 n_i g_c^{ij}\partial_j u\in H^{-1/2}(\partial B_2).
\]
We also let $\Lambda_0$ denote the ``free" DtN operator on $\partial B_2$, namely the DtN map associated with the Helmholtz equation (\ref{eq:Helm1}) with $g=I$ and $q=1$ in $B_2$. Next, we shall construct $(g_c, q_c)$ to make 
\begin{equation}\label{eq:e}
\Lambda_c=\Lambda_0.
\end{equation}
To that end, we let
\[
F(x)=\left(1+\frac{|x|}{2}\right)\frac{x}{|x|},\quad x\in B_2\backslash\{0\}
\]
It is verified that $F$ maps $B_2\backslash\{0\}$ to $B_2\backslash\overline{B_1}$. Let
\[
g_c(x)=F_* I=\frac{|x|-1}{|x|}\Pi(x)+\frac{|x|}{|x|-1}(I-\Pi(x)),\quad x\in B_2\backslash\overline{B}_1,
\]
and
\[
q_c(x)=F_* 1=\frac{4(|x|-1)}{|x|},\quad x\in B_2\backslash\overline{B}_1,
\]
where $\Pi(x):\mathbb{R}^2\rightarrow\mathbb{R}^2$ is the projection to the radial direction, defined by
\[
\Pi(x)\xi=\left(\xi\cdot\frac{x}{|x|}\right)\frac{x}{|x|},
\]
i.e., $\Pi(x)$ is represented by the symmetric matrix $|x|^{-2} xx^T$. It can be seen that $g_c$ possesses both degenerate and blow-up singularities at the cloaking interface $\partial B_1$. Hence, one need to deal with the singular Helmholtz equation (\ref{eq:Helm1})--(\ref{eq:bc}). It is natural to consider physically meaningful solutions with finite energy. To that end, one could introduce the weighted Sobolev norm (cf. \cite{LiuZho})
\begin{equation}\label{eq:weighted norm}
\|\psi\|_{g,q}=\left|  \int_{B_2}\left( g^{ij}\partial_j\psi \partial_i\psi+\omega^2 q \psi^2 \right) dx   \right|^{1/2}.
\end{equation}
It is directly verified that for $\psi\in \mathscr{E}(B_2)$
\[
\|\psi\|_{g,q}<\infty\ \ \mbox{if{f}}\ \ \frac{\partial\psi}{\partial\theta}\bigg|_{\partial B_1}=0,
\]
where $\mathscr{E}(B_2)$ denotes the linear space of smooth functions in $B_2$, and the standard polar coordinate $(x_1,x_2)\rightarrow (r\cos\theta, r\sin\theta)$ in $\mathbb{R}^2$ has been utilized. Hence, we set
\[
\mathscr{T}^\infty(B_2):=\left\{ \psi\in\mathscr{E}(B_2); \frac{\partial\psi}{\partial\theta}\bigg|_{\partial B_1}=0 \right\},
\]
and
\[
\mathscr{T}_0^\infty(B_2):=\mathscr{T}^\infty(B_2)\cap\mathscr{D}(B_2),
\]
which are closed subspaces of $\mathscr{E}(B_2)$. Then, we introduce the finite energy solution space
\begin{equation}\label{eq:space}
H_{g,q}^1(B_2):=\mbox{cl}\{ \mathscr{T}^\infty(B_2); \|\cdot\|_{g,q} \},
\end{equation}
that is, the closure of the linear function space $\mathscr{T}^\infty(B_2)$ with respect to the singularly weighted 
Sobolev norm $\|\cdot\|_{g,q}$.

(\ref{eq:e}) is justified in \cite{LiuZho} by solving (\ref{eq:Helm1})--(\ref{eq:bc}) in $H_{g,q}^1(B_2)$. Particularly, it is shown that $w=u|_{B_1}\in H^1(B_1)$ is a solution to
\begin{equation}\label{eq:interior}
\begin{cases}
& \nabla\cdot(g_a\nabla w)+\omega^2 q_a w=f_a\quad \mbox{in\ \ $B_1$},\\
& w|_{\partial B_1}=\mbox{const},\\
&\displaystyle{\int_{\partial B_1}\sum_{i,j\leq 1} n_i g_a^{ij}\partial_j w\ ds=0.}
\end{cases}
\end{equation}
In order to show the solvability of (\ref{eq:interior}) by using Fredholm theory, it is necessary to consider the following eigenvalue problem
\begin{equation}\label{eq:interior eigenvalue}
\begin{cases}
& \nabla\cdot(g_a\nabla w)+\omega^2 q_a w=0\quad \mbox{in\ \ $B_1$},\\
& w|_{\partial B_1}=\mbox{const},\\
&\displaystyle{\int_{\partial B_1}\sum_{i,j\leq 1} n_i g_a^{ij}\partial_j w\ ds=0.}
\end{cases}
\end{equation}
We let $\mathscr{R}(g_a,q_a)$ denote the set of solutions to (\ref{eq:interior eigenvalue}). Clearly, (\ref{eq:interior}) is uniquely solvable if{f} $f_a\in\mathscr{R}(g_a,q_a)^\perp$. If (\ref{eq:interior eigenvalue}) possess nontrivial solutions, then one cannot cloak a source $f\in \mathscr{R}(g_a,q_a)$ since in such a case, (\ref{eq:interior}) has no finitely energy solution, and such a radiating source would break the cloaking. Next, we note the close connection between the (\ref{eq:interior eigenvalue}) and the generalized Schiffer's conjecture (\ref{eq:general schiffer 2}). According to Definition~\ref{def:1}, it is readily seen that if the target medium parameters $(B_1,g_a,q_a)$ fails to have the Pompeiu property, then (\ref{eq:interior eigenvalue}) has nontrivial solutions. That is, one would encounter the interior resonance problem for the invisibility cloaking discussed above. But we would also like to note that the system (\ref{eq:interior eigenvalue}) arising from invisibility cloaking is still formally determined, whereas (\ref{eq:general schiffer 2}) of the generalized Schiffer's conjecture is over-determined. The connection discussed above between the singular cloaking problem and Schiffer's conjecture is natural when one is curious in whether and how such spherical cloaking design is generalized to other shaped cloaked domain.

In order to avoid the singular structure of the perfect cloaking, it is natural to incorporate regularization into the context and consider the corresponding approximate cloaking. Let
\[
F_\varepsilon(x)=\left( \frac{2(1-\varepsilon)}{2-\varepsilon}+\frac{|x|}{2-\varepsilon} \right)\frac{x}{|x|},\quad \varepsilon\in\mathbb{R}_+,
\]
which maps $B_2\backslash \overline{B}_\varepsilon \rightarrow B_2\backslash\overline{B}_1$ and let
\begin{equation}\label{eq:nonsingular}
g_c=(F_\varepsilon)_* I\qquad\mbox{and}\qquad q_c=(F_\varepsilon)_* 1.
\end{equation}
Now, we consider the nonsingular cloaking problem (\ref{eq:Helm1})--(\ref{eq:bc}) with $(g_c, q_c)$ given in (\ref{eq:nonsingular}). Let $u_\varepsilon$ denote the corresponding solution in this case. The limiting behavior of $u_\varepsilon$ as $\varepsilon\rightarrow 0^+$ was considered in \cite{LasZho,Ngu}. In order to justify the near-cloak of the regularized construction, the following eigenvalue problem arose from the corresponding study,
\begin{equation}\label{eq:int2}
\begin{cases}
\Delta w=0 \quad &\mbox{in\ $\mathbb{R}^2\backslash B_1$},\\
\nabla\cdot(g_a\nabla w)+\omega^2 q_a w=0  &\mbox{in $B_1$},\\
w^+=w^- \ \ &\mbox{on $\partial B_1$},\\
\displaystyle{\frac{\partial w^+}{\partial n}=\sum_{i,j\leq 2} n_i g_a^{ij}\partial_j w^-}\ & \mbox{on $\partial B_1$}.
\end{cases}
\end{equation}
It is shown in \cite{Ngu} that if (\ref{eq:int2}) has only trivial solution, then $u_\varepsilon|_{\mathbb{R}^2\backslash B_1}$ converges to the ``free-space" solution, whereas $u_\varepsilon|_{B_1}$ converges to $u$ implied in 
 \begin{equation}\label{eq:int3}
\begin{cases}
\Delta w=0 \quad &\mbox{in\ $\mathbb{R}^2\backslash B_1$},\\
\nabla\cdot(g_a\nabla w)+\omega^2 q_a w=f_a  &\mbox{in $B_1$},\\
w^+=w^- \ \ &\mbox{on $\partial B_1$},\\
\displaystyle{\frac{\partial w^+}{\partial n}=\sum_{i,j\leq 2} n_i g_a^{ij}\partial_j w^-}\ & \mbox{on $\partial B_1$}.
\end{cases}
\end{equation}
Otherwise, if (\ref{eq:int2}) has nontrivial solutions forming the set $W$, then it is shown that if $f_a\in W$, then the cloaking fails as $\varepsilon\rightarrow 0^+$. Now, let's take a more careful look at the eigenvalue problem (\ref{eq:int2}). By letting $z=x+i y$ and using the conformal mapping $F: z\rightarrow 1/z$, we set $v=w^+\circ F$ with $w^+=w|_{\mathbb{R}^2\backslash \overline{B}_1}$. One can verify directly that $w|_{B_1}\in H^1(B_1)$ and $v\in H^1(B_1)$ satisfy the following generalized interior transmission eigenvalue problem
\begin{equation}\label{eq:ite2}
\begin{cases}
\Delta v=0\quad &\mbox{in\ $B_1$},\\
\nabla\cdot(g_a\nabla w)+\omega^2q_a w=0\ \ &\mbox{in\ $B_1$},\\
v-w=0\quad &\mbox{on\ $\partial B_1$},\\
\displaystyle{\frac{\partial v}{\partial n}+\sum_{i,j\leq 2} n_i g_a^{ij}\partial_j w=0}\ &\mbox{on\ $\partial B_1$}.
\end{cases}
\end{equation}

From our earlier discussions in this section, one can see that the invisibility cloaking construction is very unstable, especially when there is interior resonance problem. Hence, in order to overcome this instability, cloaking schemes by incorporating some damping mechanism through the adding of a lossy layer between the cloaked and cloaking regions have been introduced and studied recently, see \cite{KOVW,LiLiuSun,LiuSun}. But it is interesting to further note that for cloaking design with a lossy layer, one may encounter eigenvalues in the complex plan, i.e., poles. In fact, it is numerically observed in \cite{LiuZhou2} for the cloaking of electromagnetic waves, the cloak effect will be deteriorated when the frequency is close to the poles.

\section*{Acknowledgement}

The work is supported by NSF grant, DMS 1207784. The author would like to thank the anonymous referee for many insightful and constructive comments.

\end{document}